\documentclass[11pt]{article}

\usepackage{a4wide}

\usepackage{amsmath}
\usepackage{amssymb}
\usepackage{amsfonts}

\date{}

\newtheorem{proposition}{Proposition}[section]
\newtheorem{theorem}[proposition]{Theorem}
\newtheorem{lemma}[proposition]{Lemma}

\newtheorem{definition}[proposition]{Definition}
\newtheorem{corollary}[proposition]{Corollary}
\newcommand{\smallbox}{{\vrule height3pt width3pt depth0pt}}
\newcommand{\qed}{\hfill \smallbox}

\newcommand{\Q}{\bigskip\par\noindent}

\def\Char{{\rm char}}

\def\mR{{\mathbb R}}
\def\mZ{{\mathbb Z}}

\def\fdeg{{\rm fdeg}}
\def\jij{_{j\in\JJ}}

\def\aut{{\rm Aut}}

\def\JJ{{\cal J}}

\def\FF{{\cal F}}

\def\EEul0{{\cal E}^{\lambda, 0}}

\def\lcc{{[[}}
\def\rc{{]]}}

\def\00{{\bf 0}}

\def\gr{{\rm gr}}

\def\aut{{\rm Aut}}

\def\nFM0{{\nu }_{F,M_0}}
\def\nFN0{{\nu }_{F,N_0}}
\def\nGN0{{\nu }_{G,N_0}}

\def\N0{ {\bf N}_0 }

\def\t{\otimes}

\def\Xpm{X^{\pm }}

\def\s{\sigma}
\def\siz{_{i\in\mZ}}

\def\l1{{\lambda}_1}

\def\m{{\bf m}}
\def\a{\alpha}
\def\a0{ {\alpha }_0}
\def\a1{ {\alpha }_1}

\def\l{\lambda}

\def\nFGM0{{\nu }_{F,G,M_0}}

\def\nFN0{{\nu}_{F,N_0}}

\def\sm{{\sigma}^m}

\def\sm1{{\sigma}^{-1}}

\def\smtp1{{\sigma}^{-t+1}}

\def\S1{S^{-1}}

\def\Xpm1{X^{\pm 1}_1}

\def\sPM1{{\sigma }^{\pm 1}}
\def\sMP1{{\sigma }^{\mp 1 }}

\def\d{\delta}

\def\Ytm1{Y^{t-1}}
\def\Yim1{Y^{i-1}}

\begin{document}

\title{Quadratic and cubic invariants of \\unipotent affine automorphisms
}
\author{V V Bavula\footnote{This research was done while the first
author held a Royal Society/NATO Fellowship
 at the University of Edinburgh.} ~and T H Lenagan}
\date{}
\maketitle


\begin{abstract}{\small Let $K$ be an  arbitrary field of
characteristic zero, $P_n:= K[ x_1, \ldots , x_n]$ be a polynomial
algebra, and $P_{n, x_1}:= K[x_1^{-1},  x_1, \ldots , x_n]$, for $n\geq
2$. Let $\s' \in {\rm Aut}_K(P_n)$ be given by
$$
x_1\mapsto x_1-1, \quad
x_1\mapsto x_2+x_1,\quad  \ldots ,\quad x_n\mapsto x_n+x_{n-1}.$$
It is proved that the algebra of invariants,
$F_n':= P_n^{\s'}$, is
a polynomial algebra in $n-1$ variables which is generated by
$[\frac{n}{2}]$  quadratic and $[\frac{n-1}{2}]$  cubic
(free) generators that are given  explicitly.

Let $\s \in {\rm Aut}_K(P_n)$ be given by
$$
x_1\mapsto x_1, \quad x_1\mapsto x_2+x_1, \quad \ldots ,\quad
x_n\mapsto x_n+x_{n-1}.$$
It is well-known that the algebra of invariants, $F_n:= P_n^\s$, is
finitely generated (Theorem of Weitzenb\"ock, \cite{Weitz},
1932), has  transcendence degree $n-1$, and that one can give an explicit
transcendence basis in which the elements have
degrees $1, 2, 3, \ldots , n-1$. However, it is an  old open problem
to find  explicit generators for $F_n$. We find
an explicit vector space
basis for the  quadratic invariants, and prove
that the algebra of invariants $P_{n, x_1}^\s$ is a polynomial
algebra over $K[x_1, x_1^{-1}]$ in $n-2$ variables which is
generated by $[\frac{n-1}{2}]$  quadratic and
$[\frac{n-2}{2}]$  cubic (free) generators that are given
explicitly.

 The coefficients of these  quadratic and cubic
invariants throw light on the `unpredictable combinatorics' of
invariants of affine automorphisms and of $SL_2$-invariants.

 \Q {\em
 }

}
\end{abstract}


\section{Introduction}

\Q Throughout the paper, $K$ denotes an {\em arbitrary}
field of characteristic zero. Let
$P_n= P^{[n]}=K[x] : = K[x_1, \dots, x_n]$ be a polynomial ring
in $n$ variables over $K$. First, we consider the $K$-algebra
automorphism $\s$ of $P_n$ given by
$$
\s : \, x_1\mapsto x_1-1, \quad x_1\mapsto
x_2+x_1, \quad \ldots , \quad x_n\mapsto x_n+x_{n-1}.
$$
This automorphism can be written
in matrix form as $\s (x) = J_n(1)x-e_1$, where $J_n(1)=
E+\sum_{i=1}^{n-1} E_{i+1, i}$ is the $n\times n$ lower triangular
Jordan  matrix ($E$ is the identity matrix and $E_{ij}$ are the
matrix units), $x= (x_1, \ldots , x_n)^t$, and $e_1=(1,0, \ldots ,
0)^t$.  It is well-known that the algebra of invariants, $P_n^\s$,
is a polynomial algebra in $n-1$ variables and that the generators can be
chosen to have degrees $2,3, \ldots , n$. (Briefly, $\s$ can be
presented as $e^\d := \sum_{i\geq 0}\frac{\d^i}{i!}$ where $\d\in
{\rm Der}_K(P_n)$ is a {\em locally nilpotent} derivation for
which there exists an element $x\in P_n$ such that $\d (x)=1$,
then $P_n^\s = P_n^\d := {\rm ker} (\d )$ and the result is old
and well-known for $\d$.) A theorem of Weitzenb\"ock \cite{Weitz}
states that the algebra of invariants $P_n^{\mathbb{G}_a}$ is
finitely generated for every linear action of the additive
(algebraic) group ${\mathbb{G}_a}$ of the field $K$ (see also
\cite{Seshadri61}, \cite{NagataL14}, and also
\cite{PopovLemSesh}). The same result is true for the algebra of
invariants $P_n^\d$ where $\d$ is a linear derivation of $P_n$; that is,
$\d (x) = Ax$ where $A$ is an $n\times n$ matrix over $K$. It is
an old open problem to find {\em explicit} generators for the
algebras $P_n^{\mathbb{G}_a}$ and $P_n^\d$. Several cases for
small $n$ are considered in \cite{NowockiLN94}.

We summarise the main results of the paper below; full proofs are given later.

The proof of the first theorem is `direct'; that is, it does not use
a reduction to the case of $\d$.

\begin{theorem}\label{Inty2}
Let $\s (x)= J_n(1) x-e_1\in {\rm Aut}_K(P_n)$, for $n\geq 2$. The
algebra of invariants $P_n^\s$ is a polynomial algebra $K[y_2,
\ldots , y_n]$ in $n-1$ variables given by
$$ y_{i+1} =\sum_{j=1}^i\phi_{-i+j}x_{j+1} + i \s^{-1}
(\phi_{-i-1}), \;\; {\rm for}\;\;
i=1, \ldots , n-1, $$ where $\phi_0:=1$ and
$\phi_{-i}:= \frac{x_1(x_1-1)\cdots (x_1-i+1)}{i!}$, for $i\geq 1$.
(Note that $\deg (y_{i+1})=i+1$.) \qed
\end{theorem}

The polynomial algebra $P_n= K[x]=\cup_{i\geq 0} K[x]_{\leq i}$ is
a filtered algebra by using the {\em total} degree of variables;
so that
$K[x]_{\leq i}:= \sum_{\deg (p)\leq i} Kp$. The {\em integer part}
of $r\in \mR$ is denoted by  $[r]:= \max\{m\in\mZ\mid m\leq r\}$.
The next theorem gives an explicit basis for the quadratic invariants
of the automorphism $\s$.
%
%
\begin{theorem}\label{usubk}
Let $\s(x) = J_n(1)x - e_1$, for $n\geq 2$, and suppose that $K$
is a field of characteristic zero.  Then the elements $u_0 =1$, and
\begin{equation*}
\begin{split}
u_k &= x_k^2 + \sum_{i=1}^{k-1}\sum_{j=k}^{2k-i}
    \,\lambda_{i,j}^k x_ix_j
+ \sum_{i=k}^{2k}\,\mu_i^kx_i,
\end{split}
\end{equation*}
where
\[\lambda_{i,j}^k = (-1)^{k-i} \left\{ {k-i\choose j-k} + {k-i-1\choose
j-k-1}\right\} \]
and
\[
 \mu_i^k = (-1)^{k-1} \left\{
{k\choose i-k} + {k-1\choose i-k-1}\right\} ,\] for $k = 1, \dots,
m:=[n/2]$, form a basis of the vector space $K[x]^\s \cap
K[x]_{\leq 2}$. In particular,
$\dim_K ( K[x]^\s \cap K[x]_{\leq 2} ) = m+1$
and $K[x]^\s \cap K[x]_{\leq 1} = K$. Each of the coefficients
$\lambda_{i,j}^{k}$ and $\mu_{i}^k$ is nonzero. \qed
\end{theorem}

\Q {\bf Remark}. ~~ In particular, $u_1= x_1^2 + x_1 + 2x_2$ and
$u_2 = x_2^2 -x_1(x_2 +2x_3) -x_2 -3x_3 -2x_4$. Note that $y_3 =
x_3 + x_1x_2 + \frac{x_1^3 - x_1}{3}$. Consider the {\em cubic}
$\s$-invariant  polynomial 
\begin{equation}\label{Intv1}
v_1 := 3y_3 = x_1^3 + 3x_1x_2  -x_1 +3x_3\in K[x_1, x_2, x_3].
\end{equation}

\begin{theorem}\label{tvsubk}
Let $\s(x) = J_n(1)x - e_1$, for $n\geq 5$, and suppose that $K$
is a field of characteristic zero.  Then, for $k=2, \dots, \mu
:=[(n-1)/2]$, the following
polynomials, $v_k$, belong to \newline
$K[x]^\s\cap K[x]_{\leq 3}$:
%

%
\begin{equation}\label{vsubn}
v_k = x_1 u_k + x_kx_{k+1} +
\sum_{i=1}^{k-1}\sum_{j=k+1}^{2k-i+1}\, \alpha_{i,j}^{k}x_ix_j +
\sum_{i=k+1}^{2k+1}\, \beta_i^k x_i,
\end{equation}
where
\begin{equation*}
\begin{split}
\alpha_{i,j}^{k}
    &= (-1)^{k-i}\,\left\{
2{k-i-1\choose j-k-1} +3{k-i-1\choose j-k-2} \right\} +
(k-i-1)\lambda_{i,j-1}^{k}\\
    &= (-1)^{k-i}\left\{ 2{k-i-1\choose j-k-1} +
3{k-i-1\choose j-k-2} + (k-i-1)\left[ {k-i\choose j-k-1} + {
k-i-1\choose j-k-2}\right]\right\},
\end{split}
\end{equation*}
and
\begin{equation*}
\beta_j^k = \alpha_{1,j-1}^{k} + \alpha_{1,j}^{k} +\mu_{j-1}^{k},
\quad {\rm for~~}j= k+1, \dots, 2k+1,
\end{equation*}
where $u_k, \lambda_{i,j}^{k}$ and $\mu_{i}^{k}$ are as defined in
Theorem~\ref{usubk}.  Note that each of  the coefficients
$\alpha_{i,j}^{k}$ and $\beta_i^k$
 is nonzero. \qed
\end{theorem}

\Q {\bf Remark}.~~  In particular, 
\begin{equation}\label{aki2k}
\alpha_{i, 2k+1-i}^k = (-1)^{k-i}(1+2(k-i))\quad {\rm and}\quad
\alpha_{i, k+1}^k = (-1)^{k-i}(k-i+1),
\end{equation}
and
\begin{equation}\label{v2exp}
 v_2 = x_1u_2 + x_2x_3
-2x_1x_3 -3x_1x_4 -3x_3 -8x_4 -5x_5.
\end{equation}

The quadratic and cubic invariants obtained in the previous two theorems
provide a generating set for the algebra of invariants, as the next theorem
shows.

\begin{theorem}\label{ctvsub}
Let $\s (x) = J_n(1)x-e_1$, for $n\geq 2$. Set  $m:= [\frac{n}{2}]$ and
$\mu :=[\frac{n-1}{2}]$. Then
\begin{enumerate}
\item $P_n^\s = K[u_1, \ldots , u_m, v_1, \ldots , v_\mu]$ is
 polynomial ring in $n-1 \; (= m+\mu )$ variables.
 \item $P_n= P_n^\s [x_1]$.
\end{enumerate}
\end{theorem}

\noindent
{\it Proof}. For each $k\geq 1$, we have 
\begin{equation}\label{Intlas}
u_k = (-1)^{k-1}2x_{2k} + \dots  \qquad {\rm and}\qquad v_k =
(-1)^{k-1}(1+2k)x_{2k+1} +\dots,
\end{equation}
where the three dots denote terms from $P_{2k-1}$ and $P_{2k}$
respectively. These imply  that
$P_n= K[u_1, \ldots , u_m, v_1, \ldots
, v_\mu][x_1]$.
It then follows that
$P_n^\s = K[u_1, \ldots , u_m, v_1, \ldots , v_\mu]$ and
$P_n= P_n^\s [x_1]$,
since $ K[u_1, \ldots , u_m, v_1, \ldots ,
v_\mu]\subseteq P_n^\s$, $\s (x_1)=x_1-1$ and ${\rm char}(K)=0$.
\qed

\Q

Now, consider the $K$-automorphism $\s (x)= J_{n+1}(1) x$ of the
polynomial algebra $P_{n+1}:= K[x_1, \ldots , x_{n+1}]$:
$$
\s :  x_1\mapsto x_1, \quad x_2\mapsto x_2+x_1, \quad \ldots , \quad
x_{n+1} \mapsto x_{n+1}+x_n.
$$
The algebra of invariants $F:= P_{n+1}^\s=\oplus_{i\geq 0}F_i$ is
a positively graded subalgebra of the polynomial algebra $P_{n+1}
= \oplus_{i\geq 0} P_{n+1, i}$ (the natural grading) where $F_i:=
F\cap P_{n+1, i}$. Let $P_{n+1, x_1}= K[x_1^{-1} , x_1, x_2,
\ldots , x_{n+1}]$ be the localization of $P_{n+1}$ at the powers
of the element $x_1$. Then $P_{n+1, x_1} = K[ x_1^\pm , z_1,
\ldots , z_n] = Q[x_1^\pm ]$ where $Q:= K[ z_1, \ldots , z_n]$ is
a polynomial algebra in the $n$ variables $z_i:= -
\frac{x_{i+1}}{x_1}$, $i=1, \ldots , n$. We denote by the same
letter $\s $ the unique extension of the automorphism $\s $ to
$P_{n+1, x_1}$. Then $\s (Q)=Q$ and $\s (z) = J_n(1)z-e_1$; that is,
$$ \s (z_1) = z_1-1, \quad \s (z_2) = z_2+z_1, \quad\ldots , \quad
\s (z_n) =
z_n +z_{n-1}.$$ Define polynomials $p_k$ and $q_k$
in $P_{n+1}$ as follows:

\begin{equation}\label{Intpk}
p_k := x_1^2u_k(z) = x_{k+1}^2 +
\sum_{i=1}^{k-1}\sum_{j=k}^{2k-i}\,
\lambda_{i,j}^{k}x_{i+1}x_{j+1} - x_1\sum_{i=k}^{2k}\, \mu_i^{k}
x_{i+1}, \;\; {\rm for}\;\; k\geq 1,
\end{equation}
while
\begin{equation}\label{Intq1}
q_1:= x_1^3 -x_2^3+3x_1x_2x_3+x_1^2x_2-3x_1^2x_4,
\end{equation}
and
\begin{equation}\label{Intqk}
q_k := x_1^3v_k(z)= -x_2p_k + x_1x_{k+1}x_{k+2} +
x_1\sum_{i=1}^{k-1}\sum_{j=k+1}^{2k-i+1}\,\alpha_{i,j}^{k}x_{i+1}x_{j+1}
- x_1^2\sum_{i=k+1}^{2k+1}\, \beta_i^kx_{i+1},
\end{equation}
for $k\geq 2$.

\begin{theorem}\label{Intpq}
Let $\s (x) = J_{n+1}(1)x$, for $n\geq 2$. Set $m:= [\frac{n}{2}]$ and
$\mu :=[\frac{n-1}{2}]$. Then the set of elements of $P_{n+1}^\s$:
$$ x_1, p_1, \ldots , p_m, q_1, \ldots , q_\mu$$
is a transcendence basis for the algebra $P_{n+1}^\s$, with $\deg
(p_i)=2$ and $\deg (q_j)=3$. Further,
$$ P_{n+1, x_1}^\s = K[x_1, x_1^{-1}][p_1, \ldots , p_m, q_1, \ldots ,
q_\mu] .$$
\end{theorem}

\noindent{\it Proof}. This follows directly from Theorem \ref{ctvsub}.
{\smallbox}

\begin{corollary}\label{Intbqd}
Let $\s (x) = J_{n+1}(1)x$, for $n\geq 2$, and set
$m:= [\frac{n}{2}]$. Then
$x_1^2,  p_1, \ldots , p_m$ is a $K$-basis of the vector space
of quadratic invariants.
\end{corollary}

\noindent
{\it Proof}. This follows from Theorem \ref{usubk} and Corollary
\ref{ngeq3}(3). {\smallbox}


\section{$\sigma$-exponentials}

\Q
Let $K$ be a field of characteristic zero, and let $\s$ denote the
affine automorphism of $K[x]$ such that $\s(x) = x-1$.  Our aim is to
choose a basis for $K[x]$ as a $K$-vector space that facilitates
calculations involving $\s$.  The idea is to exploit the fact that
$1-\s$ is a $\s$-derivation, and to choose the basis with this in mind.
Accordingly, we define
%
%
\begin{eqnarray}\label{phis}
\phi_0 := 1, \quad \phi_i := \phi_i(x) = \frac{x(x+1)\cdots (x
+i -1)}{i!} = \frac{x\s^{-1}(x) \cdots \s^{-i+1}(x)}{i!}, \;\; i\geq 1
\end{eqnarray}
and
%
%
\begin{eqnarray}\label{phi-}
\phi_0 := 1, \quad \phi_{-i} := \phi_{-i}(x)
= \frac{x(x-1)\cdots (x
-i +1)}{i!} = \frac{x\s(x) \cdots \s^{i-1}(x)}{i!}, \;\; i\geq 1.
\end{eqnarray}
Each of the two sets $\{\phi_i\}$ and $\{\phi_{-i}\}$ forms a $K$-basis of
 $K[x]$.

\Q
Note that $(1-\s)\phi_i = \phi_{i-1}$ and  $(1-\s)\phi_{-i} =
\s(\phi_{-i+1})$, for $i\geq 1$, while $(1-\s)\phi_0 = 0$.  Note also
that  $\phi_{-i}(-x) = (-1)^i\phi_i(x)$, and that $\sigma^{i-1}(\phi_i)
= \phi_{-i}$, for $i\geq 1$.

\Q The choice of bases and the action of $1-\s$ suggests that we should
construct exponential functions, twisted by $\s$.  In order to do this,
we extend the automorphism $\s$ to an automorphism of the
the power series ring $K[x][[\Theta]]$
by defining $\s(\Theta) =\Theta$.

\Q
Now, define
%
%
\begin{eqnarray}\label{E}
E= E(x) := \sum_{i=0}^{\infty}\, \phi_{-i}\Theta^i =  1 +
\sum_{i=1}^{\infty}\, x\s(x) \cdots
\s^{i-1}(x)\frac{\Theta^i}{i!}
\end{eqnarray}
and
%
%
\begin{eqnarray}\label{E-}\lefteqn{
E_-= E(-x) := \sum_{i=0}^{\infty}\, \phi_{-i}(-x)\Theta^i = }\\
&&\sum_{i=0}^{\infty}\, (-1)^i\phi_i\Theta^i =
1 + \sum_{i=1}^{\infty}\, (-1)^ix\s^{-1}(x) \cdots
\s^{-i+1}(x)\frac{\Theta^i}{i!}.
\end{eqnarray}
The following identities are easily established by direct computation:
 \[
(1-\s)E= \Theta\s(E),\qquad  \text{ and }\qquad (1-\s)E_- = -\Theta
E_-
.\]

\begin{lemma} $
E(x)^{-1} = E(-x) $ in $K[x]\lcc \Theta\rc$.
\end{lemma}
\Q
{\em Proof}.~~ Set $E=E(x) $ and $E_- = E(-x)$.
By applying the $\s$-derivation $(1-\s)$ to the product $E_-E$,
and,  by using the identities  $(1-\s)E= \Theta\s(E)$ and $(1-\s)E_- =
 -\Theta E_-$, we obtain
%
%
\begin{eqnarray}\label{EE-}
(1-\s)(E_-E) &=& (1-\s)E_-\cdot E +\s(E_-)(1-\s)E \\
&=& -\Theta(1-\s) (E_-E).
\end{eqnarray}
It follows that $ (1-\s)(E_-E) \in \cap_{n=1}^{\infty} \Theta^n
K[x][[\Theta]] =0$
 and so
$\s(E_-E) = E_-E$.  Hence, $E_-E \in k[x][[\Theta]]^{\s} = K[[\Theta]]$,
and we may write
 $E_-E = 1 + \sum_{i=1}^{\infty}\, \lambda_i\Theta^i $, with each
$\lambda_i \in K$.

\Q
By setting $x=0$ in the previous equality,  we get
\[
1=1\cdot 1 = E_-(0)E(0) = 1 + \sum \lambda_i\Theta^i,
\]
and it follows that each $\lambda_i =0$; so that $E_-E =1$ in
$K[x]\lcc \Theta\rc$. {\smallbox}

\Q Consider the
$K$-automorphisms $\s_i$ of the polynomial ring $K[x_1,x_2]$ in
two variables, defined by $ \s_i(x_j):= x_j -\delta_{ij}$,  for
$i,j = 1,2$, where $\delta_{ij}$ is the Kronecker delta symbol.
The automorphisms $\s_i$ extend uniquely to automorphisms of the
algebra $K[x_1,x_2]\lcc\Theta\rc$, by setting $\s_i(\Theta)
=\Theta$.
%
%
\begin{lemma}\label{powerlaw}
$ E(x_1)E(x_2) = E(x_1 + x_2) $ in $K[x]\lcc\Theta\rc$.
\end{lemma}

\Q {\em Proof}.~~  It suffices to show that the product $P:=
E_-(x_1)E_-(x_2) E(x_1 + x_2) $ is equal to one.   Note that the
identity $(1-\s_i)(E(x_1 +x_2)) = \Theta\s_i(E(x_1+x_2))$ holds,
since $\s_i(x_1 +x_2) = (x_1 +x_2) -1$. Hence, by  using the same
argument as in the proof of (\ref{EE-}), one easily obtains $
(1-\s_i)P = (-\Theta)^n(1-\s_i)P$,  for all $n\geq 1$, and
$i=1,2$. Hence,
\[
P = 1+\sum_{i=1}^{\infty}\, \lambda_i\Theta^i \in \cap_{i=1}^{2}
\ker(1-\s_i) = K\lcc \Theta\rc,
\]
so that each $\lambda_i \in K$.
Now, set $x_1=x_2=0$ in the previous equality, to obtain
\[
1 = E_-(0)E_-(0)E(0) = 1 +\sum_{i=1}^{\infty}\, \lambda_i\Theta^i.
\]
Hence, all $\lambda_i =0$; and so $P=1$, as required.
{\smallbox}

\Q
The following useful identity now follows immediately.
%
%
\begin{lemma}\label{=0}
1.~~ For all $k\geq 1$,
\[
(-1)^k\sum_{\substack{i+j=k\\i\geq 0,\; j\geq 0}}\,
(-1)^j\phi_i\phi_{-j} = \sum_{\substack{i+j=k\\i\geq 0,\; j\geq
0}}\, (-1)^i\phi_i\phi_{-j} =0.
\]

\Q 2.~~ For all $n, k \geq 1$,
\[
(-1)^k\sum_{\substack{i+j=k\\n \geq i\geq 0,\; n\geq j\geq 0}}\,
(-1)^j\phi_i\phi_{-j} = \sum_{\substack{i+j=k\\n \geq i\geq 0,\;
n\geq j\geq 0}}\, (-1)^i\phi_i\phi_{-j} =0.
\]

\end{lemma}

\Q
{\em Proof}.~ 1.~~ This follows immediately from the equality
\[
1=E_-E  = (\sum_{i\geq 0} (-1)^i\phi_{i} \Theta^i)(\sum_{j\geq 0}
\phi_{-j} \Theta^j) = \sum_{k\geq 0 }\left( \sum_{i+j=k}\,
(-1)^i\phi_i\phi_{-j}\right)\Theta^k. \hspace{5ex}
\]

\Q
2.~~ Follows immediately from the equality above by working modulo
$\theta^n$.  {\smallbox}

\Q In order to study the Jordan blocks occurring in the canonical
form of an affine automorphism, we need to
 consider specializing the above results to the case that $\Theta$
is the nilpotent \mbox{$(n-1)\times (n-1)$} matrix
\[
\left(
\begin{array}{cccccc}
0 & 0 & 0 & 0 &\cdots & 0\\
1 & 0 & 0& 0 &\cdots & 0\\
0 & 1 & 0 & 0&& 0\\
  &   &\ddots    & \ddots & & \\
  &   &   & \ddots &\ddots & \\
0 & 0 & 0  & \cdots & 1 & 0
\end{array}
\right)
.\]
Note that $\Theta^{n-1}= 0$, but $\Theta^{n-2} \neq 0$.
\Q
Consider the matrix
\begin{align}
\Lambda &=  \sum_{i=0}^{n-2}(-1)^i\phi_i\Theta^i \\[3ex]
    & =
\left(
\begin{array}{ccccccc}
1 & 0 & 0 &  \cdots &\cdots & 0 & 0\\
-\phi_1 & 1& 0 &   \cdots &\cdots & 0 & 0\\
\phi_2  & -\phi_1  & 1 &\ddots    &  & & \\
-\phi_3  &  \phi_2  & -\phi_1 &  & \ddots & & \\
\vdots & \vdots  &   & \cdots & &&\\
\vdots & \vdots &&&&&\\
\vdots & \vdots &&&&&\\
(-1)^{n-2}\phi_{n-2} & (-1)^{n-3}\phi_{n-3} & \cdots \cdots & &&-\phi_1 &1
\end{array}
\right) \in {\rm SL}_{n-1}(K[x])
.\end{align}
The above analysis reveals that
\[
\Lambda^{-1} = \sum_{i=0}^{n-2}\, \phi_{-i} \Theta^i.
\]

\Q Set $\Phi :=(-\phi_2, \phi_3, \dots, (-1)^i\phi_{i+1}, \dots,
(-1)^{n-1}\phi_n)^t \in K[x]^{n-1}$. \qed

%
%
\begin{lemma}\label{eta}
\[
\Lambda^{-1}\Phi =
\left( \begin{array}{c}
-\s^{-1}(\phi_{-2})\\
\vdots\\
-i\s^{-1}(\phi_{-i-1})\\
\vdots\\
-(n-1)\s^{-1}(\phi_{-n})
\end{array}\right)
, \]
where the $i$th entry is displayed.
\end{lemma}
\Q
{\em Proof}.~~
Set

\[
 \left(\begin{array}{c} \eta_1\\
\vdots \\ \vdots \\\eta_{n-1}
\end{array}\right) :=
\Lambda^{-1}\Phi =
\left(
\begin{array}{ccccccc}
1 & 0 & 0 &  \cdots &\cdots & 0 & 0\\
\phi_{-1} & 1& 0 &   \cdots &\cdots & 0 & 0\\
\phi_{-2}  & \phi_{-1}  & 1 &\ddots    &  & & \\
\phi_{-3}  &  \phi_{-2}  & \phi_{-1} &  & \ddots & & \\
\vdots & \vdots  &   & \cdots & &&\\
\phi_{-i+1} &  &&&&&\\
\vdots & \vdots &&&&&\\
\phi_{-n+2} & \phi_{-n+3} & \cdots \cdots & &&\phi_{-1} &1
\end{array}
\right)
\left(\begin{array}{c}
-\phi_2\\
\phi_3\\
\vdots\\
(-1)^{i}\phi_{i+1}\\
\vdots\\
(-1)^{n-1}\phi_n
\end{array}\right)
.\]
Then

\begin{eqnarray*}
\eta_i &=&
    \sum_{j=1}^{i}\, \phi_{-i+j}(-1)^{j}\phi_{j+1} =
-\sum_{l=2}^{i+1}\, (-1)^l\phi_{l}\phi_{-i-1+l} \\
    &=&
    -\sum_{l=0}^{i+1}\, (-1)^l\phi_{l}\phi_{-i-1+l} +\phi_{-i-1}
-\phi_1\phi_{-i} \\
    &=&
    0 + \phi_{-i}\left(\frac{x-i}{i+1} -x\right) =
-i\frac{(x+1)\phi_{-i}}{i+1} = -i\s^{-1}(\phi_{-i-1}),
\end{eqnarray*}
as claimed.  Note that the $0$ in the  above calculation
arises by applying Lemma~\ref{=0}.
 {\smallbox}


\section{The invariant polynomials $u_i$ and $v_j$}

Let $K$ be a field of characteristic zero, and let $\s$ be the
affine automorphism of $K[x] := K[x_1, \ldots , x_n]$, for $n\geq
2$,  defined by the rule
$$ \s (x_1)= x_1-1, \quad \s (x_2) = x_2+x_1, \quad\ldots , \quad\s (x_n) =
x_n+x_{n-1}. $$
In matrix form,
\[
\s(x) = J_n (1)x -e_1,
\]
where $J_n(1)= E+\sum_{i=1}^{n-1} E_{i+1, i}$ is the $n\times n $
lower triangular Jordan matrix ($E$ is the identity matrix and
$E_{ij}$ are the matrix units).  Observe that $\s (K[x_1, \ldots
,x_m] ) = K[x_1, \ldots ,x_m]$, for each $i\geq 1$. In particular,
$\s(K[x_1]) = K[x_1]$, and $\s(x_1) = x_1 -1$.

\Q Recall, from the
previous section, that the set of polynomials $\{\phi_i := \phi_i
(x_1)\}$ defined by
%
%
\begin{eqnarray}\label{phi}
\phi_0 := 1, \quad \phi_i := \phi_i(x_1) = \frac{x_1(x_1+1)\cdots (x_1
+i -1)}{i!}, \;\; i\geq 1
\end{eqnarray}
is a $K$-basis of $K[x_1]$, and that $(1-\s)\phi_i = \phi_{i-1}$, for all
$i\geq  1$, while $(1 -\s)\phi_0 =0$.
\Q
The matrix $\Theta := J_{n-1}(1) - I$ is a nilpotent matrix with
$\Theta^{n-1} =0$, but $\Theta^{n-2} \neq 0$.  As in the previous
section, we consider the
matrix
\begin{align*}\Lambda &=  \sum_{i=0}^{n-2}(-1)^i\phi_i\Theta^i\\
    &=
\left(
\begin{array}{ccccccc}
1 & 0 & 0 &  \cdots &\cdots & 0 & 0\\
-\phi_1 & 1& 0 &   \cdots &\cdots & 0 & 0\\
\phi_2  & -\phi_1  & 1 &\ddots    &  & & \\
-\phi_3  &  \phi_2  & -\phi_1 &  & \ddots & & \\
\vdots & \vdots  &   & \cdots & &&\\
\vdots & \vdots &&&&&\\
\vdots & \vdots &&&&&\\
(-1)^{n-2}\phi_{n-2} & (-1)^{n-3}\phi_{n-3} & \cdots \cdots & &&-\phi_1
&1
\end{array}
\right)
 \in {\rm SL}_{n-1}(K[x_1])
.\end{align*}
\Q
Set
\[
x' = (x_2, \ldots,  x_{i+1}, \ldots, x_{n})^t \quad {\rm  and}\quad
\Phi = (-\phi_2, \phi_3, \dots, (-1)^i\phi_{i+1}, \dots,
(-1)^{n-1}\phi_n)^t
.\]

\Q
Define $y= (y_2,\ldots, y_{n})^t
 \in K[x]^{n-1}$ by the linear equation
$x' = \Lambda y +\Phi$; so that $y = \Lambda^{-1} (x' - \Phi)$.
In more detail, we have
%
%
\begin{eqnarray}\label{xi}
x_{i+1}  &=& \sum_{j=1}^{i}\, (-1)^{i-j}\phi_{i-j}y_{j+1} +
(-1)^{i}\phi_{i+1}
\end{eqnarray}
and
%
%
\begin{eqnarray}\label{yi}
y_{i+1} &=& \sum_{j=1}^{i}\, \phi_{-i+j}x_{j+1} +
i\s^{-1}(\phi_{-i-1}),
\end{eqnarray}
for $i= 1, \ldots , n-1$, by Lemma~\ref{eta}.
We extend the action of $\s$ to the $(n-1)\times (n-1)$ matrix ring
$M_{n-1}(K[x])$ and to the column space $K[x]^{n-1}$ in the obvious way
(that is, elementwise).

\Q
The following proposition contains the claim of Theorem~\ref{Inty2}.
%
%

\begin{proposition}\label{1fixed}
Let $\s(x) = J_n(1)x -e_1$, for $n\geq 2$, and suppose that $\Char
(K)=0$.  Then, the fixed ring $K[x]^{\s}$ is equal to the polynomial ring
$K[y_2, \ldots, y_{n}]$ in the $n-1$ variables  defined by (\ref{yi}).
 Further, $K[x] = K[x_1]\t K[x]^\s = K[x_1, y_2,  \ldots, y_{n}]$.
\end{proposition}

\Q
{\em Proof}.\,  Note that the subalgebra $K[y]$ of $K[x]$, generated  by
$y_2, \ldots , y_{n}$, is isomorphic to a polynomial ring in $n-1$
variables, by (\ref{yi}).

\Q
Observe that
\[
(1-\s)x' = -\Theta x' - (\phi_1, 0, \dots , 0)^t, \qquad (1-\s)\Phi =
-\Theta \Phi - (\phi_1, 0, \dots , 0)^t
\]
and
\[
(1-\s)\Lambda^{-1} = \Theta\s(\Lambda^{-1}).
\]
By applying the $\s$-derivation $(1-\s)$ to the equation
$y= \Lambda^{-1}(x' - \Phi)$, we obtain
\begin{eqnarray*}
(1-\s)y &=& (1-\s)\Lambda^{-1}\cdot (x' -\Phi) +
    \s(\Lambda^{-1})\cdot (1-\s)(x'-\Phi)\\
    &=& \Theta \s(\Lambda^{-1})(x' -\Phi) -
\s(\Lambda^{-1})\Theta(x' - \Phi) =0,
\end{eqnarray*}
since $\Theta\s(\Lambda^{-1}) = \s(\Lambda^{-1})\Theta$.

\Q
Thus, $(1-\s)y=0$,
 and so $\s(y_i) = y_i$, for each $i=2, \dots, n$.
Hence, $k[y] := k[y_2, \dots, y_{n}] \subseteq k[x]^{\s}$, and
$k[x] = k[y] \otimes k[x_1]$.

\Q
Let $f= \sum_{i=0}^{s} f_i\phi_i \in K[x]^{\s}$, where each $f_i \in
K[y]$.  Then,
\[
0= (1-\s)f = f_1\phi_0 +\cdots + f_s\phi_{s-1},
\]
and it follows that $f_i = 0$, for all $i\geq 1$.  Thus, $f=f_0 \in
K[y]$, and $K[x]^{\s} \subseteq K[y]$, as required.  {\smallbox}

\Q
Let $K$ be a commutative ring and let $Z = I\times \mZ$ be a subset of
$\mZ^2$, where $I= [a, a+1, \dots, b]$, for some $a<b$.  Suppose that
$\lambda, \mu:Z \rightarrow K$ are functions such that $(i,j)
\rightarrow \lambda_{i,j}$ and $(i,j) \rightarrow \mu_{i,j}$ and such
that the relation
\[
\lambda_{i,j} = \delta(\lambda_{i+1, j-1} + \lambda_{i+1,j}) +
\mu_{i,j-1},
\]
holds for all $(i,j) \in Z$ and for some $\delta \in K$, then we write
$\mu\stackrel{\delta}{\rightsquigarrow} \lambda$

%
%
\begin{lemma}\label{Fun}
Let $K$ be a commutative ring and suppose that two
 functions $\lambda, \mu :Z
\rightarrow K$ given by $(i,j) \mapsto \lambda_{i,j}$  and $(i,j)
\mapsto \mu_{i,j}$ satisfy the
relation
$\lambda_{i,j} = \delta(\lambda_{i+1,j-1} + \lambda_{i+1, j}) +
\mu_{i,j-1}$,
for all $(i,j)\in Z$ and
for some $\delta \in K$.  Then,

\Q
1.
%
%
\begin{eqnarray}\label{Fun1}
\lambda_{i,j} = \delta^c \sum_{d=0}^{c}\, {c\choose d}\lambda_{i+c,j-d}
+ \sum_{c'=0}^{c-1}\delta^{c'} \sum_{d'=0}^{c'}\,
 {c'\choose d'}\mu_{i+c', j-1-d'},
\end{eqnarray}
for each integer $c\geq 0$, such that $i+c \in [a,b]$.
\Q
In particular, when $\mu =0$, we have
%
%
\begin{eqnarray}\label{Fun2}
\lambda_{i,j} = \delta^c \sum_{d=0}^{c}\, {c\choose d}\lambda_{i+c,j-d},
\end{eqnarray}
for each integer $c\geq 0$, such that $i+c \in [a,b]$. \newline
\noindent 2.
If, in addition, the function $\mu$ satisfies the relation
$
\mu_{i,j} = \gamma(\mu_{i+1,j-1} + \mu_{i+1,j}),
$
for all $(i,j)\in Z$ and for some unit $\gamma \in K$, then
%
%
\begin{eqnarray}\label{Fun3}
\lambda_{i,j} = \delta^c \sum_{d=0}^{c}\, {c\choose d}\lambda_{i+c,j-d}
+ \left(1 + \frac{\delta}{\gamma} + \left(\frac{\delta}{\gamma}\right)^2
+ \cdots + \left(\frac{\delta}{\gamma}\right)^{c-1}\right)\mu_{i,j-1}.
\end{eqnarray}
\end{lemma}

\Q
{\em Proof}.~~ 1.~~
We use induction on $c$.  The base cases of the induction $c=0,1$ are
obvious. Suppose that $c\geq 2$ and that the result holds for $c-1$.
Denote by $\Omega = \Omega_c$ the second sum in (\ref{Fun1}).
Then, by induction,
\begin{equation*}
\begin{split}
\lambda_{i,j} &= \delta^{c-1} \sum_{d=0}^{c-1}\,
{c-1 \choose d}\lambda_{i+c-1,j-d}
+ \sum_{c'=0}^{c-2}\delta^{c'} \sum_{d'=0}^{c'}\,
 {c'\choose d'}\mu_{i+c', j-1-d'}\\
    &= \delta^c \sum_{d=0}^{c-1}\, {c-1\choose d}(
\lambda_{i+c,j-d-1} + \lambda_{i+c,j-d})
+ \sum_{c'=0}^{c-1}\delta^{c'} \sum_{d'=0}^{c'}\,
 {c'\choose d'}\mu_{i+c', j-1-d'} \\
    &= \delta^c \sum_{d=0}^{c}\,\left\{
 {c-1\choose d} + {c-1\choose d-1}\right\}\lambda_{i+c,j-d}
+ \Omega\\
    &= \delta^c \sum_{d=0}^{c}\, {c\choose d}\lambda_{i+c,j-d}
+\Omega .
\end{split}
\end{equation*}
as required.

\Q 2.~~By (\ref{Fun2}), $ \mu_{i,j-1} =
\gamma^{c'}\sum_{d=0}^{c'}\,{c'\choose d'}\,u_{i+c', j-1-d'} $,
for each integer $c'\geq 0$.  The element $\gamma$ is a unit, so
\[
\Omega_c = \left(1 + \frac{\delta}{\gamma} +
\left(\frac{\delta}{\gamma}\right)^2
+ \cdots + \left(\frac{\delta}{\gamma}\right)^{c-1}\right)\mu_{i,j-1},
\]
and the equation (\ref{Fun3}) follows.
{\smallbox}

\Q
{\bf Remark}.~~   We are setting  ${a\choose b} =
0$, for each pair $a,b\in \mZ$ that does not satisfy $0\leq b\leq a$.

\begin{corollary}
Let $K$ be a commutative ring.  Suppose that the functions $\lambda,
\lambda^1, \dots, \lambda^n:Z\rightarrow K$ satisfy
\[
0\stackrel{\delta_{n}}{\rightsquigarrow}\lambda^{n}
\stackrel{\delta_{n-1}}{\rightsquigarrow}\lambda^{n-1}
\stackrel{\delta_{n-2}}{\rightsquigarrow} \dots
\stackrel{\delta_{1}}{\rightsquigarrow}\lambda^{1}
\stackrel{\delta_{0}}{\rightsquigarrow}\lambda^{0} \equiv \lambda,
\]
where $\delta_1, \dots, \delta_n$ are units in $K$. \newline
Then \newline
1.
\[
\lambda_{i,j} = \delta_0^c\sum_{d=0}^{c}\, {c\choose d}\lambda_{i+c,j-d}
+ \sum_{k=1}^{n}\, (-1)^{k-1}\lambda_{i,j-k}^{k} \left( \sum_{c_1 =
0}^{c-1}\sum_{c_2 = 0}^{c_1-1}\dots \sum_{c_k =
0}^{c_{k-1}-1}
\prod_{l=1}^k\left(\frac{\delta_{l-1}}{\delta_l}\right)^{c_l}\right),
\]
for each integer $c\geq 0$ such that $i+c\in [a,b]$.

\Q
2.~~ In particular, when $\delta_0 = \delta_1 = \dots = \delta_n$, we
have
%
%
\begin{align}\label{Fun4}
\lambda_{i,j} &= \delta_0^c\sum_{d=0}^{c}\, {c\choose d}\lambda_{i+c,j-d}
+ \sum_{k=1}^{\min(n,c)}\, (-1)^{k-1} \phi_k(c-k+1)\lambda_{i,j-k}^k\\
    &= \delta_0^c\sum_{d=0}^{c}\, {c\choose d}\lambda_{i+c,j-d} +
\sum_{k=1}^{n}\,(-1)^{k-1} {c\choose k}\lambda_{i,j-k}^k.
\end{align}
Moreover,
%
%
\begin{align}\label{Fun5}
\lambda_{i,j} &= \delta_0^c\sum_{d=0}^{c}\, {c\choose d} \left(
\sum_{k=0}^{\min(n,c)}\, \phi_k(c)\lambda_{i+c,j-k-d}^k \right)\\
    &= \delta_0^c\sum_{k=0}^{\min(n,c)}\, {c+k-1\choose k}
\left\{ \sum_{d=0}^{c}\, {c\choose d}\lambda_{i+c,j-k-d}^k
\right\} .
\end{align}
\end{corollary}

\Q {\em Proof}.~~ 1.~~We use induction on $n$.   The base case
$n=1$ was proved in (\ref{Fun3}).  Suppose that $n\geq 2$, and
that the result holds for the case $n-1$.  By induction, we have
\[
\lambda_{i,j-1}^1 = \delta_1^{c_1}
\sum_{d'=0}^{c_1}\, {c_1\choose d'}\lambda_{i+c_1,j-1-d'}^1
+ \sum_{k=1}^{n-1}\, (-1)^{k-1}\lambda_{i,j-1-k}^{k+1}
\left( \sum_{c_2 =
0}^{c_1-1}\dots \sum_{c_{k+1} =
0}^{c_{k}-1}
\prod_{l=1}^k\left(\frac{\delta_{l}}{\delta_{l+1}}\right)^{c_{l+1}}
\right),
\]
for each integer $c_1\geq 0$ such that $i+c_1\in [a,b]$. Combining
the above equality with (\ref{Fun1}) in the case $\lambda^{1}$;
that is, with
$\lambda^1\stackrel{\delta_{0}}{\rightsquigarrow}\lambda$,
\[
\lambda_{i,j} = \delta_0^c \sum_{d=0}^{c}\, {c\choose d}\lambda_{i+c,j-d}
+ \sum_{c_1=0}^{c-1}\delta_0^{c_{1}} \sum_{d'=0}^{c_1}\,
 {c_1\choose d'}\lambda_{i+c_1, j-1-d'}^1,
\]
we obtain
\begin{align*}
\lambda_{i,j}
    &=
\delta_0^c \sum_{d=0}^{c}\, {c\choose d}\lambda_{i+c,j-d}\\
 & \hspace{5ex} + \sum_{c_1 =0}^{c-1}
    \left( \frac{\delta_{0}}{\delta_1} \right)^{c_1}
    \left\{
\lambda_{i,j-1}^{1} -  \sum_{k=1}^{n-1}\,
(-1)^{k-1}\lambda_{i,j-1-k}^{k+1}
    \left(
\sum_{c_2
=
0}^{c_1-1}\dots \sum_{c_{k+1} =
0}^{c_{k}-1}\prod_{l=1}^k
    \left(
\frac{\delta_{l}}{\delta_{l+1}}
        \right)^{c_{l+1}}
        \right)
        \right\}\\
    &=
 \delta_0^c\sum_{d=0}^{c}\, {c\choose d}\lambda_{i+c,j-d}
+ \sum_{k=1}^{n}\, (-1)^{k-1}\lambda_{i,j-k}^{k}
    \left(
\sum_{c_1 =
0}^{c-1}\sum_{c_2 = 0}^{c_1-1}\dots \sum_{c_k =
0}^{c_{k-1}-1}\prod_{l=1}^k
    \left(
\frac{\delta_{l-1}}{\delta_l}
    \right)^{c_l}
    \right).
    \end{align*}

\Q
2.~~ If $\delta_0 = \dots = \delta_n$, we will prove by induction on $k$
that
%
%
\begin{equation}\label{multsum}
I_k:=\sum_{c_1 =
0}^{c-1}\sum_{c_2 = 0}^{c_1-1}\dots \sum_{c_k =
0}^{c_{k-1}-1}
\prod_{l=1}^k\left(\frac{\delta_{l-1}}{\delta_l}\right)^{c_l}
= \begin{cases} \phi_k(c-k+1), & \text{if $c\geq k$}\\
        0, & \text{if $c<k$}. \end{cases}
\end{equation}
Obviously, $I_k = 0$, whenever $c<k$, so we assume that $c\geq k$.
Now,
\[
I_k = \sum_{c_1=0}^{c-1}\, \phi_{k-1}(c_1 -k +2) = \phi_{k-1}(1)
+\phi_{k-1}(2) + \dots + \phi_{k-1}(c-k+1) = \phi_{k}(c-k+1),
\]
since $0, -1, -2, \dots, -(k-2)$ are roots of the polynomial
$\phi_{k-1}$ and $c\geq k$. Since $\sigma^{k-1}(\phi_k) = \phi_{-k}$,
for all $k\geq 1$, we see that
\[
\phi_{k}(c-k+1) = \left( \sigma^{k-1}(\phi_k) \right) (c) =
\phi_{-k}(c) = \frac{c(c-1) \dots  (c-k+1)}{k!},
\]
for $c\geq k$.  Hence, $I_k = {c\choose k}$, since we are setting
${a\choose b} =0$, for each pair $a,b\in\mZ$ that does not satisfy
$0\leq b \leq a$.  Thus, the formula (\ref{Fun4}) follows.

\Q
Since
\[
\phi_k(c) = \frac{c(c+1) \dots (c+k-1)}{k!} = {c+k-1\choose k},
\]
for $k\geq 0$, the second equality in (\ref{Fun5}) follows from the
first; so, it remains to prove the first equality in (\ref{Fun5}).  We
use induction on $n$.  The case $n=0$ is evident, see (\ref{Fun2}), so
we suppose that $n\geq 1$.
By (\ref{Fun4}),
\[
\lambda_{i,j} = \delta_0^c\sum_{d=0}^{c}\, {c\choose
d}\lambda_{i+c,j-d}
- \sum_{k=1}^{\min(n,c)}\, (-1)^{k} \phi_{-k}(c)\lambda_{i,j-k}^k
.\]
By induction on $n$, setting $m:= \min(n,c)$, we have
\begin{align*}
\lambda_{i,j}
    &=
\delta_0^c\sum_{d=0}^{c}\, {c\choose
d}\lambda_{i+c,j-d} - \sum_{k=1}^{\min(n,c)}\, (-1)^k\phi_{-k}(c)
    \left\{
\delta_0^c\sum_{d=0}^{c}\, {c\choose
d}
    \left(
\sum_{l=0}^{\min(n-k,c)}\,\phi_l(c) \lambda_{i+c,j-(k+l)-d}^{k+l}
    \right)
    \right\}\\
    &=
\delta_0^c\sum_{d=0}^{c}\, {c\choose
d}
    \left\{
    \lambda_{i+c, j-d} - \sum_{k=1}^{\min(n,c)}\,
\sum_{l=0}^{\min(n-k,c)}\, (-1)^k \phi_{-k}(c)\phi_{l}(c)
\lambda_{i+c,j-(k+l)-d}^{k+l}
    \right\}\\
    &=
\delta_0^c\sum_{d=0}^{c}\, {c\choose
d}
        \left\{
        \lambda_{i+c, j-d} - \sum_{s=1}^{m}
    \left(
\sum_{\substack{k+l =s\\m\geq k \geq 1, m\geq l \geq 0}}\,
(-1)^k\phi_{-k}(c)\phi_l(c)
    \right)\lambda_{i+c, j-s-d}^{s}
    \right\}\\
    &=
\delta_0^c\sum_{d=0}^{c}\, {c\choose
d}
        \left\{
        \lambda_{i+c, j-d} + \sum_{s=1}^{\min(n,c)}\,
\phi_s(c)\lambda_{i+c,j-s-d}^s
    \right\},
\end{align*}
    as required,
since
\[
\sum_{\substack{k+l =s\\m\geq k \geq 1, m\geq l \geq 0}}\, (-1)^k
\phi_{-k}(c)\phi_l(c)
    =
-\phi_s(c) + \sum_{\substack{k+l =s\\m\geq k \geq 0, m\geq l \geq 0}}\,
(-1)^k \phi_{-k}(c)\phi_l(c)
        = -\phi_s(c) +0 = -\phi_s(c),
\]
by Lemma~\ref{=0}.   {\smallbox}

\Q
We are now in a position to prove Theorem~\ref{usubk} in which we find a
basis for the quadratic invariants.

\Q {\bf  Proof of Theorem \ref{usubk}}. Any  element of the set
$K[x]_{\leq 2}$ which has constant term equal to zero  can  be
written as a sum
\[
u= \sum_{j=1}^{n} \lambda_j x_j^2 + \sum_{i=1}^{n-1}x_i\left( \sum_{j=
i+1}^{n} \lambda_{ij}x_j \right) +\sum_{j=1}^{n} \mu_jx_j.
\]
The element $u$ is uniquely determined by the upper triangular $n\times
n$ matrix $\Lambda = (\lambda_{ij}) \in M_n (K)$ and the vector $(\mu_1,
\dots, \mu_n) \in K^n$.  For the sake of convenience, we will set
$\lambda_i = \lambda_{ii}$.

\Q
Observe that $u \in K[x]^\s $ if and only if
$\partial (u) =0$, where $\partial =
1-\s$.  In order to calculate $\partial(u)$, we perform elementary
computations using the fact that $\partial$ is a $\s$-derivation, and the
following facts: $\partial(x_1) = 1$ and $\partial (x_1^2)
 = 2x_1 -1$; while
$\partial(x_i) = -x_{i-1}$ and $\partial(x_i^2)
 = -2x_{i-1}x_i - x_{i-1}^2$, for
$i= 2, \dots, n$.  We obtain
%
%
\begin{equation}\label{d(u)}
\begin{split}
\partial(u) &= -\sum_{j=1}^{n-1}(\lambda_{j,j+1} + \lambda_{j+1})x_j^2 \\
    &~~~-\sum_{i=2}^{n-2} x_i\left\{ ( \lambda_{i,i+2} +
\lambda_{i+1, i+2} + 2\lambda_{i+1})x_{i+1}
+ \sum_{j=i+2}^{n-1}(\lambda_{i+1,j} + \lambda_{i, j+1}
+\lambda_{i+1,j+1})x_j
 + \lambda_{i+1,n}x_n \right\}\\
    &~~-x_1\left\{ (\lambda_{2,3} + \lambda_{1,3}+2 \lambda_{2})x_2 +
\sum_{j=3}^{n-1}(\lambda_{2,j} + \lambda_{1,j+1} +\lambda_{2,j+1})x_j +
\lambda_{2,n}x_n\right\}\\
    &~~~ -2\lambda_{n}x_{n-1}x_n + (2\lambda_{1} +\lambda_{1,2}
-\mu_2)x_1 \\
    &~~~ +\sum_{j=3}^{n}(\lambda_{1,j-1} + \lambda_{1,j}
-\mu_j)x_{j-1} +\lambda_{1,n}x_n + (\mu_1 - \lambda_{1}).
\end{split}
\end{equation}
Thus, $\partial(u) =0$ if and only if each of the coefficients in the
expression above are zero.  This gives the system of linear equations
below (see (\ref{mu}), (\ref{c1}),  (\ref{c2}),  (\ref{c3}), (\ref{c4})
below).

\Q
We can immediately see from the coefficients that the entries in the
last column of the matrix $\Lambda$ must all be zero for a solution to
$\partial(u) = 0$.  Also, the linear terms are specified by the last few
coefficients, viz:
%
%
\begin{equation}\label{mu}
\mu_1 = \lambda_1, \quad \mu_2 = 2\lambda_1 + \lambda_{1,2}, \quad
\mu_j= \lambda_{1, j-1} + \lambda_{1,j}, \quad j = 3, \dots, n.
\end{equation}
The remaining equations can be separated into four classes:
%
%
\begin{equation}\label{c1}
\lambda_{2,3} + \lambda_{1,3} + 2\lambda_{2} =0,
\end{equation}
%
%
\begin{equation}\label{c2}
\lambda_{j,j+1} + \lambda_{j+1} = 0,\quad j = 1, \dots, n-1,
\end{equation}
%
%
\begin{equation}\label{c3}
\lambda_{i+1,j} + \lambda_{i,j+1} + \lambda_{i+1,j+1} = 0,\quad i=1,
\dots, n-2, \quad j= i+2, \dots, n-1,
\end{equation}
%
%
\begin{equation}\label{c4}
 \lambda_{i, i+2} + \lambda_{i+1, i+2} + 2\lambda_{i+1} =0, \quad i=2,
\dots, n-2.
\end{equation}

\Q
Obviously, the elements $u_i$ are linearly independent; so, it suffices
to prove that an element $u\in K[x]^\s \cap K[x]_{\leq 2}$,
with zero constant term,
is a linear combination of the elements $u_i$.  In order to do this, we
will use induction on $n\geq 2$.  In the case that $n=2$, we see that
the element $u_1 = x_1^2 + x_1 +2x_2$ is the unique solution (up to
non-zero scalar multiple) of the system $\partial(u) =0$.  The same is true
for $n=3$, since $\lambda_{2} = \lambda_{12}=0$, by using (\ref{c1}) and
(\ref{c2}).

\Q
Thus, we may assume that $n\geq 4$, and that the result is true for all
$n'$ strictly less than $n$.

\Q
The last column of the matrix $(\lambda_{ij})$ is zero.  By using
(\ref{c3}) with $i=1, \dots, n-2$ and $j = n-1$, we see that
$\lambda_{i,n-1} =0$, for $i=2, \dots, n-2$.  Since $\lambda_{n-1}=
-\lambda_{n-2, n-1} = 0$, by (\ref{c2}) with $j=n-2$, it follows that
$\lambda_{i,n-1}$ =0, for all $i>1$.  By using similar arguments, it
follows that all of the elements of the matrix $(\lambda_{ij})$ lying
below and on the anti-diagonal are zero; that is,
%
%
\begin{equation}\label{R}
\lambda_{i,j} =0, \quad i+j \geq n+1.
\end{equation}
By passing from $u$ to a suitable linear combination of the form
 $
u+\sum_{i=1}^{m'} \, \alpha_iu_i, $ where $m' = [(n-1)/2]$ and
$\alpha_i \in K$, we may assume that
%
%
\begin{equation}\label{L}
\lambda_{1} = \cdots = \lambda_{m'} =0.
\end{equation}

\Q In more detail, we will solve the system assuming that the
conditions above hold, as a result we will have the polynomials
$u_i$, and this justifies our assumption.  By (\ref{c2}), we have
\[
\lambda_{j,j+1} =0, \quad j = 1, \dots, m'-1,
\]
and
\[
\lambda_{j,j+2} = 0,\quad j =1\dots, m'-2,
\]
by (\ref{c1}) and (\ref{c4}).
Now, by (\ref{c3}), we obtain
%
%
\begin{equation}\label{LL}
\lambda_{i,j} =0, \quad i, j = 1, \dots, m' .
\end{equation}
In the case that $n$ is even, it is enough to prove that the element
$u_m$ is the unique solution (up to a nonzero scalar multiple)
satisfying (\ref{L}), and in the case that $n$ is odd, that the only
solution satisfying (\ref{L}) is $0$.

\Q
Suppose first that $n$ is even, with $n =2m$, and $m\geq 2$ (the cases
where $n =2,3$ have been considered earlier).
Suppose that $n=4$.  By (\ref{c2}), we see that $\lambda_{1,2} =
-\lambda_{2}$, and by (\ref{c1}), we obtain $\lambda_{1,3} =
-2\lambda_{2}$; and so, $u = \lambda_2u_2$, by (\ref{mu}).

\Q
Suppose now that $n=6$.  By (\ref{c2}), we see that $\lambda_{2,3} =
-\lambda_{3}$, and by (\ref{c1}), we obtain $\lambda_{1,3} =
\lambda_{3}$.  Now, $\lambda_{2,4} =-2\lambda_{3}$, by (\ref{c4}).  By
(\ref{c3}), $\lambda_{1,4}= 3\lambda_{3}$ and $\lambda_{1,5}=
2\lambda_{3}$, hence $u = \lambda_3u_3$, by (\ref{mu}).  Finally (for
the even case), suppose that $n=2m$, with $m\geq 4$.  In this case,
$m'=m-1$, so all of the diagonal elements of the matrix $(\lambda_{ij})$,
are zero, except for $\lambda_m = \lambda_{mm}$.  By (\ref{c2}), we
obtain $\lambda_{m-1,m} = -\lambda_m$, and, by using (\ref{c4}), we get
$\lambda_{m-2,m} = \lambda_m$; then, by (\ref{c2}),
\begin{equation}
\lambda_{i,m} = (-1)^{m-i}\lambda_m, \quad i=1, \dots, m.
\end{equation}
By (\ref{c4}), $\lambda_{m-1,m+1} = -2\lambda_m$, and then, by
(\ref{c3}),
\begin{equation}
\lambda_{m-i,m+i} = (-1)^i2\lambda_m,\quad i =1, \dots, m-1.
\end{equation}
Now, all of the entries of the first $m-1$ rows satisfy (\ref{c3}), and
$\lambda_{m-1,m} = -\lambda_m, \lambda_{m-1,m+1} = -2\lambda_m$; also,
all of the other entries of the $(m-1)$-st row are zero.  If we apply
Lemma~\ref{Fun} to  the entries of the first $m-1$ rows, and put $\delta
= -1$, we obtain
%
%
\begin{equation}\label{lij}
\begin{split}
\lambda_{ij} &= (-1)^{m-1-i} \sum_{k=0}^{m-1-i}\,
{m-1-i\choose k}\lambda_{m-1,j-k}\\
    &= (-1)^{m-1-i}\lambda_m\left\{ - {m-1-i\choose j-m}
-2{m-1-i\choose j-m-1} \right\}\\
    &= (-1)^{m-i}\lambda_m\left\{  {m-i\choose j-m}
+{m-1-i\choose j-m-1} \right\},
\end{split}
\end{equation}
for $i = 1, \dots, m-1,\; j = m, \dots, 2m-i$.

\Q
Thus,
%
%
\begin{equation}\label{muj}
\begin{split}
\mu_j &= (-1)^{m-1}\lambda_m\left\{ {m-1\choose j-m-1} + {m-2\choose
j-m-2}   +
{m-1\choose j-m} + {m-2\choose j-m-1}   \right\}\\
    &= (-1)^{m-1}\lambda_m\left\{ {m\choose j-m}+ {m-1\choose
j-m-1}   \right\},
\end{split}
\end{equation}
by (\ref{mu}) and (\ref{lij}).
This finishes the even case.

\Q
Now, suppose that $n$ is odd, with $n =2m+1$, for some $m\geq 2$.   In
this case, $m' =m$, and so $\lambda_{m+1} =0$, by (\ref{R}); and
$\lambda_{m,m+1} =0$, by (\ref{c2}).
If $m=2$, then $\lambda_{1,3} =0$, by (\ref{c1}), and $\lambda_{1,4}
=0$, by (\ref{c3}); so that $\lambda_{ij} =0$, for all $i,j$.  If $m\geq
3$, then $\lambda_{m-1,m+1} =0$, by (\ref{c4}).  It then follows that all
$\lambda_{ij}$ are zero, by using (\ref{R}), (\ref{L}) and (\ref{c3}).
{\smallbox}

\Q Theorem~\ref{tvsubk}, which presents the  cubic invariants, can
now be proved.

\Q
 {\bf  Proof of Theorem \ref{tvsubk}}.  The element $u_k$
of Theorem~\ref{usubk} can be written as the  sum,  $u_k = u'_k +
u''_k$, of the quadratic terms $u'_k$ and the linear terms
$u''_k$.  The element $v_k = x_1u_k + v'_k + v''_k$, where
\[
v'_k = ax_kx_{k+1} + \sum_{i=1}^{k-1}\sum_{j=k+1}^{2k-i+1}\,
\alpha_{i,j}^{k}x_ix_j
\quad {\rm and}\quad
v''_k = \sum_{i= k+1}^{2k+1}\, \beta_i^k x_i.
\]

\Q
Clearly,
\[
\partial(v_k) = u_k +
\partial(v'_k) + \partial(v''_k) = (u'_k + \partial(v'_k))
+ (u''_k + \partial(v''_k))
.\]

\Q
Thus, by using (\ref{d(u)}),
we see that  $\partial(v_k)=0$ if and only if
the coefficient $a$ and the $\alpha_{i,j}^k, \beta_i^k$
satisfy the following
system of linear equations:
\begin{align*}
-a+1 &= 0\\
-a-\alpha_{k-1,k+1}^{k} -1 &= 0\\
-a -\alpha_{k-1,k+2}^{k} -2 &= 0\\
-(\alpha_{i+1,j}^{k} + \alpha_{i,j+1}^{k} +\alpha_{i+1,j+1}^{k} )
+\lambda_{i,j}^{k} &=0,
\end{align*}
for $i= 1, \dots, k-1$ and $j=k, \dots, 2k-i$;
and
%
%
\[\label{betaj}
\beta_j^k = \alpha_{1,j-1}^{k} + \alpha_{1,j}^{k} +\mu_{j-1}^{k},
\]
for $j = k+1, \dots, 2k+1$ (note that we set
$\alpha_{1,2k+1}^{k} =0$).

\Q
Equivalently,
\[
a=1,\quad \alpha_{k-1,k+1}^k =-2, \quad \alpha_{k-1,k+2}^k =-3
\]
and
\[
\alpha_{i,j}^{k} = \delta( \alpha_{i+1,j-1}^{k} +
\alpha_{i+1,j}^{k})+\lambda_{i,j-1}^{k},
\]
for $i= 1, \dots, k-1$ and $j=k+1, \dots, 2k-i+1$, where $\delta = -1$.

\Q We know, from the proof of Theorem~\ref{usubk}, that
$\lambda_{i,j}^{k} = \delta(\lambda_{i+1,j-1}^{k} +
\lambda_{i+1,j}^{k})$, for $i= 1, \dots, k-1$ and $j=k, \dots, 2k-1$.
Thus, by using Lemma~\ref{Fun}.(2), we have
\begin{equation*}
\begin{split}
\alpha_{i,j}^{k}
    &= (-1)^{k-i-1}\left\{ {k-i-1\choose j-k-1}\alpha_{k-1,
k+1}^{k} + {k-i-1\choose j-k-2}\alpha_{k-1, k+2}^{k}\right\}
+ (k-i-1)\lambda_{i,j-1}^{k}\\[3ex]
    &=
(-1)^{k-i}\left\{  2{k-i-1\choose j-k-1} +
 3{k-i-1\choose j-k-2} +
(k-i-1)\left[ {k-i\choose j-k-1} + { k-i-1\choose j-k-2}\right]\right\}
.~~~~{\smallbox}
\end{split}
\end{equation*}


\section{$\FF$-direct sums}

\Q Let $Q = \cup_{i\in \mZ}\, Q_i$ be a $\mZ$-filtered algebra
with a filtration $\FF = \{ Q_{i}\}$.  We will always assume that
the filtration is {\bf separated}; that is, $\cap_{i\in\mZ}\,
Q_{i} =0$. Any subspace $U$ of $Q$ has an induced filtration $U =
\cup_{i\in\mZ}\, U_{i}$, where $U_i:= U\cap Q_i$.  In this case,
the associated graded space $\gr_{\FF} U = \oplus\siz\,
U_i/U_{i-1}$ is a {\em subspace} of the associated graded algebra
$\gr_{\FF} Q = \oplus\siz\,  Q_i/Q_{i-1}$ in a natural manner.

\Q
Given a separated filtration  $\FF = \{ Q_{i}\}$ of $Q$, then, for any
nonzero element $u\in Q$, there exists a unique $i\in\mZ$ such that
$u\in Q_i \backslash Q_{i-1}$.  The integer $i$ is called the {\bf
$\FF$-degree} of $u$, and is denoted by $\fdeg(u)$.

\begin{definition} Let $\{U_j, j\in \JJ\}$  be a set of subspaces of the
$\FF$-filtered algebra $Q$.  We say that the sum $\sum_{j\in\JJ}\,
U_j$ is {\bf $\FF$-direct} if $ \sum\jij\, \gr_{\FF} U_j =
\bigoplus\jij\, \gr_{\FF} U_j $ in $\gr_{\FF}Q$.
\end{definition}
The concept of $\FF$-directness is extremely useful in finding a
$K$-basis of a ring of invariants and in  proving that relations
for a ring of invariants are {\em defining} relations.  For a
separated filtration $\FF$ it follows easily that any $\FF$-direct
sum $\sum_{j\in\JJ}\, U_j$ is the direct sum,
 $\oplus\jij \,
U_j$,  of the subspaces $U_j$.

\begin{lemma}
Let $Q= \cup\siz \, Q_i$ be a filtered algebra with separated
filtration $\FF = \{Q_i\}$ and  $\sum_{j\in\JJ}\, U_j$ be an
$\FF$-direct sum of subspaces $\{U_j\}$.  Then
    \begin{enumerate}
\item if $u_j\in U_j$, then $\fdeg(\sum u_j) = \max \{
\fdeg(u_j)\}$.

\item  $(\sum\jij\, U_j) \cap Q_i = \sum\jij\, U_j \cap Q_i.$
    \end{enumerate}
\end{lemma}

\Q {\em Proof}.~~ 1. This is evident.

\Q
2. Denote by $L$ and $R$ the left and right hand side
vector spaces in the
equality that
we are trying to establish.  Clearly, $L\supseteq R$.  If $u =
\sum u_j \in (\sum U_j) \cap Q_i$, for some $u_j \in U_j$, then each
$u_j \in U_j \cap Q_i$, by statement 1, so $R\subseteq L$.
{\smallbox}

\Q

A $K$-basis $\{ U_i, i\in J\}$ of the filtered algebra $Q=\cup_{i\in
\mathbb{Z}} Q_i$ is called an $\FF$-{\em basis} if the sum of
$1$-dimensional subspaces $\sum_{j\in J} Ku_j$ is $\FF$-{\em
direct}. In this case, $\{ \gr \, u_j, j\in J\}$ is a $K$-basis
for the associated graded algebra $\gr \, Q:= \oplus_{i\in
\mathbb{Z}} Q_i/ Q_{i+1}$. If, in addition, the algebra
$Q=\cup_{i\geq 0} Q_i$ is {\em positively} graded then the
converse is true: a basis $\{ u_j, j\in J\}$ of $Q$ is an
$\FF$-basis of $Q$ if and only if
$\{ \gr \, u_j, j\in J\}$ is a basis for $\gr \,
Q$; and a basis $\{ u_j, j\in J_i\}$ of $Q_i$ is an $\FF$-basis
of $Q$ if and only if
$\{ \gr \, u_j, j\in J_i\}$ is a basis for $\oplus_{\nu
=0}^iQ_\nu / Q_{\nu -1}$.

\Q Similarly, elements $\{ u_j, j\in J\}$ of $Q$ are $\FF$-{\em
independent} if the sum of $1$-dimensional subspaces $\sum_{j\in
J} Ku_j$ is $\FF$-direct. In this case, elements $\{ \gr \, u_j,
j\in J\}$ are linearly independent elements of $\gr \, Q$. The
converse is obviously true; so, the elements $\{ u_i, i\in J\}$ of
$Q$ are $\FF$-independent if and only if the elements
$\{ \gr \, u_j, j\in J\}$ are
linearly independent in $\gr \, Q$.

\section{Number of variables $\leq 5$}

\Q
Let $K$ be a field of characteristic zero.  The polynomial ring, $P\equiv
P^{[n+1]} \equiv K[x] \equiv K[x_1, \dots, x_{n+1}]$, in $n+1$ variables,
 is a positively graded $K$-algebra $P = \oplus_{i\geq 0}\, P_i$, where
$P_i \equiv K[x]_i$ consists of the homogeneous polynomials of degree
$i$, together with zero.

\Q Consider the graded automorphism $\sigma \in \aut_{\rm gr}(P)$
defined by $\sigma(x) = \JJ_{n+1}(1)x$, where $\JJ_{n+1}(1)$ is
the $(n+1)\times (n+1)$ lower triangular Jordan matrix with $1$ in
each diagonal entry; that is,
$$ \s (x_1)=x_1, \quad \s (x_2)= x_2+x_1, \quad\ldots , \quad\s (x_{n+1})=
x_{n+1}+x_n.$$

\Q We use the results in the earlier sections of the paper to give
explicit generators and defining relations for $P^{\sigma}$ for
some small values of $n$.

\Q The algebra $F= P^{\sigma}$ of invariants is a positively
graded algebra $F = \oplus_{i\geq 0}\, F_i$, where $F_i=P_i \cap
F$. In the case that $n=0$ we have that $\sigma$ is the identity
map, so that $F = K[x_1]$.  Thus, we assume that $n\geq 1$.
  The element $x_1$ is $\s$-invariant.
 Denote by $P_{x_1}$ the localization of $P$ at the powers of
$x_1$,  that is,
\[
P_{x_1} = S^{-1}K[x] = K[x_1,x_1^{-1}, x_2,  \ldots , x_{n-1},
x_n, x_{n+1}] ,\] where $S= \{ x_1^{i} \mid i\geq 0\}$. Set $z_i:=
-\frac{x_{i+1}}{x_1}$, for $i= 1, \ldots, n$,  so that
%
%
\begin{eqnarray}\label{Px1}
P_{x_1} = K[z, x_1^{\pm 1}] = K[z_1, \ldots, z_{n}, x_1^{\pm 1}] =
Q[x_1^{\pm 1}],
\end{eqnarray}
where $Q = Q^{[n]} = K[z] = K[z_1, \dots, z_n]$ is the polynomial
ring in $n$ variables.  The algebra $Q = \oplus_{i\geq 0}\, Q_i$
is a positively graded $K$-algebra, using the degree of the
polynomials.  The filtration $\FF = \{ Q_{\leq i} := \oplus_{j\leq
i}\, Q_j\}$, for $i\geq 0$, associated with this grading,
satisfies $Q_{\leq i} = P_ix_1^{-i}$, for $i\geq 0$,  and so $Q=
\sum_{i\geq 0}\, P_ix_1^{-i}$.

\Q
Let $p(x_1, \dots, x_{n+1}) \in P$ be a homogeneous polynomial.  Then
%
%
\begin{eqnarray}\label{pxpz}
p(x_1, \dots, x_{n+1}) = (-x_1)^{\deg(p)}p(-1, -\frac{x_2}{x_1}, \dots,
-\frac{x_{n+1}}{x_1}) = (-x_1)^{\deg(p)}p(-1, z_1, \dots, z_n),
\end{eqnarray}
where $p(-1, z_1, \dots, z_n) \in Q_{\leq \deg(p)}$.

\Q Denote by the same letter $\sigma$ the {\em unique extension}
of the automorphism $\sigma $ to the localized algebra $P_{x_1}$.
Then $\sigma (Q) = Q$, and $\sigma(z) = \JJ_n(1)z -e_1$; that is,
$$ \s
(z_1)=z_1-1, \; \s (z_2)= z_2+z_1, \ldots , \s (z_n)=
z_n+z_{n-1}.$$
Theorem~\ref{Inty2} (or Proposition~\ref{1fixed})
now becomes available for use.

\subsection*{The case $n=1$}

\Q
Clearly, \[
\left( K[x_1, x_2]_{x_1} \right)^\sigma = K[z_1, x_1^{\pm 1}]^\sigma =
K[x_1^{\pm 1}],
\]
since $\sigma(z_1) = z_1 -1$ and the characteristic of $K$  is zero.
Hence,
%
%
\begin{equation}\label{fix1}
K[x_1, x_2]^\sigma = \left( K[x_1, x_2]_{x_1} \right)^\sigma \cap K[x_1,
x_2] = K[x_1].
\end{equation}
Thus, we may assume that $n\geq 2$.
 The first part of the next corollary follows
immediately from Proposition~\ref{1fixed} (statement 3 follows
from (\ref{pxpz}); statements 4 and 5 follow from the definition of
$\FF$-basis).
%
%
\begin{corollary}\label{ngeq3}
  Let $n\geq 2$.  Then
\newline
1.~~$Q^{\s} = K[y_2, \ldots,
y_{n}]$ is a polynomial ring in the $n-1$ variables $y_i$ given by
\[
y_{i+1} = y_{i+1}(z) =
 \sum_{j=1}^{i}\, \phi_{-i+j}(z_{1})z_{j+1} +
i\s^{-1}(\phi_{-i-1}(z_1)),
\]
for $i=1,\ldots, n-1$, where the  $\phi_l$ are defined in
(\ref{phi-}) and $Q = K[z_1]\otimes Q^\s$.

\Q
2.~~$y_i \equiv z_i \pmod {z_1}$, for $i=2, \dots, n$.

\Q 3.~~$F= \oplus_{i\geq 0}\, F_i$, with $F_i = x_1^iQ_{\leq
i}^\sigma$, where $Q_{\leq i}^\sigma = Q^\sigma \cap Q_{\leq i}$,
for $i\geq 0$.
 \Q
4. If $\{ b_j , j\in J_i\}$ is an $\FF$-basis for the vector space
$Q^\s_{\leq i}$ then $\{ x_1^{i-\deg_z (b_j)} (x_1^{\deg_z (b_j)}
b_j), j\in J_i\}$ is a $K$-basis for $F_i$. If $\{ b_j , j\in
J\}$ is an $\FF$-basis for the vector space $Q^\s$ then $\{
x_1^{i-\deg_z (b_j)} (x_1^{\deg_z (b_j)} b_j), j\in J\}$ is a
$K$-basis for $F_i$.

\Q 5.~~Given $g_1, \ldots , g_m\in Q^\s$ such that, for each
$i\geq 0$, $\{ g^\alpha := g_1^{\alpha_1}\cdots g_m^{\alpha_m}\, |
\, \alpha =(\alpha_1, \ldots , \alpha_m)\in J_i\subseteq
\mathbb{N}^m\}$ is an $\FF$-basis for $Q^\s_{\leq i}$ then $\{
x_1^{i-\sum_{j=1}^m\alpha_j\deg_z (g_j)} \prod_{j=1}^m(x_1^{\deg_z
(g_j)} g_j)^{\alpha_j}, j\in J_i\}$ is a $K$-basis for $F_i$;
and if $J_1\subseteq J_2\subseteq \cdots $ then $$\{
x_1^{i-\sum_{j=1}^m\alpha_j\deg_z (g_j)} \prod_{j=1}^m(x_1^{
\deg_z (g_j)} g_j)^{\alpha_j}\, | \,  j\in \cup_{i\geq 1}J_i, \;
i-\sum_{j=1}^m\alpha_j\deg_z (g_j)\geq 0 \}$$ is a $K$-basis for
$F$.  ~~~{\smallbox}
\end{corollary}

%
%
\begin{corollary}\label{u=uy}
1.~~Let $u = u(z_1, \dots, z_n)\in Q^\s$.  Then
\[
u(z_1, \dots, z_n) = u(0, y_2, \dots, y_n).
\]

\Q 2.~~ In particular, for the elements $u_k$ and $v_k$ in
Theorems~\ref{usubk} and~\ref{tvsubk},
\begin{align*}
u_k &= y_k^2 + \sum_{i=2}^{k-1}\sum_{j=k}^{2k-i}\, \lambda_{i,j}^k
y_iy_j + \sum_{i=k}^{2k}\,\mu_i^ky_i,\\
v_k &= y_ky_{k+1}+\sum_{i=2}^{k-1}\sum_{j=k+1}^{2k-i+1}\,
\alpha_{i,j}^{k}y_iy_j +\sum_{i=k+1}^{2k+1}\, \beta_i^ky_i.
\end{align*}
\end{corollary}

\Q
{\em Proof}.~~ Note that $Q = Q^\s\oplus Qz_1,\;\;
 Q^\s = K[y_2, \dots, y_n]$
and $y_i \equiv z_i \pmod {z_1}$, for $i=2, \dots, n$, by
Corollary~\ref{ngeq3}.  Thus, for any $u(z_1, \dots, z_n) \in Q^\s$, we
have
\[
u(z_1, z_2, \dots, z_n) = u(0, y_2, \dots , y_n) + z_1v ,\] fore
some polynomial $v\in Q$.  However,
$$ u(z_1, \ldots , z_n)- u(0, y_2, \ldots , y_n)= vz_1\in Q^\s
\cap Qz_1=0,$$ since $Q= Q^\s\oplus Qz_1$. Hence, $u(z_1, \dots,
z_n) = u(0,y_2, \dots, y_n)$.

\Q
2.~~ Evident.  {\smallbox}

\Q
The polynomials
\[
f_{i+1}:= x_1^{\deg(y_{i+1})}y_{i+1} = x_1^{i+1}y_{i+1},
\]
for $i=1, \dots, n-1$, belong to the  algebra $F$ of invariants.
Note that
\begin{align*}
f_{i+1} &= \sum_{j=1}^{i}\, (-1)^{i-j+1}x_1^j\frac{x_2(x_2+x_1)\dots
(x_2 + (i-j-1)x_1)}{(i-j)!}x_{j+2}\\
    &~~~ +
(-1)^{i+1}i\frac{(x_2-x_1)x_2(x_2+x_1)\dots
(x_2 + (i-1)x_1)}{(i+1)!}.
\end{align*}

%
%
\begin{corollary}\label{f=fn}
For each homogeneous polynomial $f(x_1, \dots, x_{n+1}) \in F = K[x_1,
\dots, x_{n+1}]^\s$, where $\s(x) = J_{n+1}(1)x$:
\[
f(x_1, \dots, x_{n+1}) = (-x_1)^{\deg(f)}f(-1, 0, f_2/x_1^2, \dots,
f_n/x_1^n).
\]
\end{corollary}

\Q
{\em Proof}.~~  This follows immediately from the equality (\ref{pxpz})
and Corollary~\ref{u=uy}.(1):
\begin{align*}
f(x_1, \dots, x_{n+1}) &= (-x_1)^{\deg(f)}f(-1,z_1, \dots, z_n)\\
        &=(-x_1)^{\deg(f)}f(-1,0, y_2, \dots, y_n)\\
        &=(-x_1)^{\deg(f)}f(-1,0, f_2/x_1^2, \dots,
f_n/x_1^n).  ~~~~~{\smallbox}
\end{align*}

\subsection*{The case $n=2$}

\Q By Corollary~\ref{ngeq3}.(1), we know that $Q^\sigma = K[u_1]$,
 where $ u_1 := 2y_2
= z_1^2 + z_1 + 2z_2$. Clearly, $\{ u_1^i, i\geq 0\}$ is an
$\FF$-basis for $Q^\s$. By Corollary~\ref{ngeq3}(5),  the algebra
%
%
\begin{equation}\label{fix2}
K[x_1, x_2,
x_3]^\sigma = K[x_1, p_1],
\end{equation}
 is the polynomial ring in the two
variables $x_1$ and $p_1 := x_1^2u_1^2 = x_2^2 - x_1(x_2 +2x_3)$.

\subsection*{The case $n=3$}
\Q
Recall that
\[
y_2 = z_2 + \frac{(z_1 + 1)z_{1}}{2}\quad {\rm and}\quad y_3 = z_3 +
 z_1z_2
+ \frac{z_1^3 - z_1}{3}
.\]
Set
%
%
\begin{align}\label{v1}
v_1 &:= 3y_3 = z_1^3 + 3z_1z_2  -z_1 +3z_3 = z_1^3 + \dots,
\end{align}
Clearly,
\begin{eqnarray*}
 v_1&=& z_1u_1-z_1^2+z_1z_2-z_1+3z_2 \\
 &=& z_1u_1-u_1+z_1z_2+5z_2.
\end{eqnarray*}
Consider the element
%
%
\begin{eqnarray}\label{omega}
\theta := v_1^2 -u_1^3 +3v_1u_1 +2u_1^2.
\end{eqnarray}
Direct computation gives
\begin{equation}
\begin{split}
\theta &= -3z_1^2 \{z_2^2 -z_1(z_2 +2z_3)\} + 9z_1(z_1 + 2z_2)z_3
-8z_2^3 +z_1z_2(5z_1+6z_2)\\
    &\quad + 9z_3^2 + 3(z_1 + 6z_2)z_3 + 8z_2^2 + 2z_1z_2\\
    &=-3z_1^2 \{z_2^2 -z_1(z_2 +2z_3)\} +\cdots .
\end{split}
\end{equation}

\Q The leading terms of the elements $u_1= z_1^2+\cdots $
 and $\theta$ are {\em algebraically
independent},  so the subalgebra $U:= K[ u_1, \theta]$ of $Q^\s =
K[u_1, v_1] $ is isomorphic to a polynomial ring in two variables,
and the elements $\{ u_1^i\theta^j\, | \, i,j\geq 0\}$ are
$\FF$-independent.

\begin{lemma}\label{CFQf}
The sum $Q^\s = U + Uv_1$
 is an $\FF$-direct sum in the filtered algebra $Q$
with the filtration $\FF = \{Q_{\leq i}\}$. In particular, $Q^\s =
U \oplus Uv_1$.
\end{lemma}

\Q
{\em Proof}.~~
    It follows from (\ref{omega}) that $Q^\s = U + Uv_1$.
Observe that the degree of the leading term of any element from
$U$ is {\em even}, and  the degree of the leading term of any
element of $Uv_1$ is {\em odd}. The result then follows.
{\smallbox}

\begin{corollary}\label{cCFQf}
For $i\geq 0$, we have
\[
Q_{\leq i}^\s = U_{\leq i} \oplus U_{\leq (i-3)}\,v_1,
\]
where $U_{\leq i} = U \cap Q_{\leq i} = \bigoplus\{ Ku_1^{i_2}
\theta^{i_4} \mid i_2, i_4\geq 0, 2i_2+4i_4 \leq i\}$, and $\{
u_1^{i_2} v_1^{i_3}\theta^{i_4} \mid  2i_2+3i_3+4i_4 \leq i;  \;
i_3=0,1; \;  i_2, i_4\geq 0\}$ is an $\FF$-basis for $Q^\s_{\leq
i}$.
\end{corollary}

\begin{theorem}
The fixed algebra $F= K[x_1, x_2,x_3,x_4]^\s$
is generated by the four elements
\[
x_1,\quad p_{1} := x_1^2u_{1} = x_2^2 -x_1(x_2 + 2x_3),\quad q_1
:= x_1^3v_1 = -x_2^3 + 3x_1x_2x_3 + x_1^2x_2-3x_1^2x_4
\]
and
\begin{multline*}
s: = x_1^4\theta  = -3x_2^2\left\{ x_3^2 - x_2(x_3+2x_4)\right\} -
9x_1x_2(x_2 + 2x_3)x_4 + 8x_1x_3^3\\ - x_1x_2x_3(5x_2 + 6x_3) +
9x_1^2x_4^2 + 3x_1^2(x_2 + 6x_3)x_4 + 8x_1^2x_3^2 + 2x_1^2x_2x_3.
\end{multline*}
of degrees $1,2,3,4$, respectively, subject to the
defining relation
%
\begin{align}\label{relation}
x_1^2s &= q_1^2 + 3x_1p_1q_1 -p_1^3 +2x_1^2p_1^2.
\end{align}
For $i\geq 0$,  $F_i=\oplus \{ K x_1^{i_1} p_1^{i_2} q_1^{i_3}
s^{i_4} \, | \, i_1+2i_2+3i_3+4i_4 = i;  \; i_3=0,1; \; i_1, i_2,
i_4 \geq 0\}$.
\end{theorem}

\Q {\em Proof}.~~ By Corollary \ref{ngeq3}, $F = \oplus_{i\geq
0}\, F_i$, where $F_i = x_1^i Q_{\leq i}^\s$. By Corollary
\ref{ngeq3}(5) and Corollary \ref{cCFQf}, (using $\FF$-directness)
$F_i=\bigoplus Kx_1^{i_1} p_1^{i_2}q_1^{i_3} s^{i_4}$ where all
$i_\nu \geq 0$,
and $i_3=0, 1$ while
$i_1+2i_2+3i_3+4i_4=i$. We
conclude that the fixed algebra $F$ is generated by the four
elements above, subject to the  relation (\ref{relation}).  The
relation is {\em the}  defining relation for $F$ (as the algebra
that satisfies this defining relation has the {\em same} basis $\{
x_1^{i_1} p_1^{i_2}q_1^{i_3} s^{i_4}\}$ as $F$).  {\smallbox}

\subsection*{The case $n=4$}

\Q Recall that (Corollary \ref{ngeq3} and Theorem \ref{usubk})
\[
y_4 = z_4 + z_1z_3 + \frac{z_1(z_1 -1)}{2!}z_2 + \frac{3(z_1
+1)z_1(z_1-1)(z_1-2)}{4!},
\]
%
%
and
\begin{equation}\label{etalead}
\begin{split}u_2 &= z_2^2 -z_1(z_2 + 2z_3) -z_2 -3z_3 -2z_4= z_2^2 -z_1(z_2 +
2z_3)+\cdots .
\end{split}
\end{equation}
By Corollary \ref{u=uy},
%
%
\begin{equation}\label{u2}
\begin{split}u_2&:= y_2^2 -y_2 -3y_3 -2y_4 = \frac{1}{4}u_1^2
-\frac{1}{2}u_1 -v_1 -2y_4.
\end{split}
\end{equation}
Consider the element
%
%
\begin{equation}\label{tilde}
\begin{split}
\tilde{\theta}  &:= \theta +3u_1 u_2 = v_1^2 -u_1^3 +3v_1u_1 + 2u_1^2
+ 3u_1u_2.
\end{split}
\end{equation}
Direct computation gives

%
%
\begin{equation}\label{tildelead}
\begin{split}
\tilde{\theta}
    &= -6z_1^2z_4 + 6z_1(z_2 -z_1)z_3 -2z_2^3 + z_1z_2(3z_2 -z_1) \\
    &\quad -6(z_1 + 2z_2)z_4 + 9z_3^2 - 6z_1z_3 + 2z_2^2 -z_1z_2\\
    &=-6z_1^2z_4 + 6z_1(z_2 -z_1)z_3 -2z_2^3 + z_1z_2(3z_2
    -z_1)+\cdots .
\end{split}
\end{equation}

\Q The leading terms of the elements $u_{1}=z_1^2+\cdots$,
 $u_2=z_2^2-z_1(z_2+2z_3)+\cdots$, and $\tilde{\theta}$  are {\em
algebraically independent};  so the algebra $Q^\s = K[y_{2},
y_{3}, y_{4}]$ contains the subalgebra $U:= K[u_1, u_2,
\tilde{\theta}]$ which is
 a polynomial ring in three variables.

\Q
For $l\geq 0$, we set
\[
U_{\leq l} := U\cap Q_{\leq l} =
 \bigoplus_{i,j, k \geq 0}\,
\{ Ku_1^iu_2^j\tilde{\theta}^k \mid 2i + 2j + 3k  \leq l\}.
\]
The set $\{ u_1^iu_2^j\tilde{\theta}^k\, | \, 2i + 2j + 3k  \leq
l\}$ is an $\FF$-basis for $U_{\leq l}$, and $\{
u_1^iu_2^j\tilde{\theta}^k\, | \,i,j,k\geq 0\}$ is an $\FF$-basis
for $U$. Now,
\begin{equation}\label{U+Uycase4}
Q^\s = K[u_1, u_2, \tilde{\theta}] +
K[u_1, u_2, \tilde{\theta}]y_{3} = U +
Uv_1.
\end{equation}
In more detail,
\begin{eqnarray*}
 Q^\s &=& K[y_2, y_3, y_4]= K[u_1, v_1, y_4]= K[u_1, v_1, u_2] \;\; ({\rm by}\; (\ref{u2})) \\
 &=& K[u_1, u_2,\tilde{\theta}]+ K[u_1, u_2,\tilde{\theta}]v_1 \;\; ({\rm by}\;
 (\ref{tilde})).
\end{eqnarray*}

\begin{lemma}\label{f7}
The sum $Q^\s = U + Uv_1$
 is an $\FF$-direct sum in the filtered algebra $Q$
with the filtration $\FF = \{Q_{\leq i}\}$. Hence,
$ \{ u_1^mu_2^j\tilde{\theta}^kv_1^l\, | \,
2m+2j+3k+3l\leq i, \; l=0,1\}$ is an $\FF$-basis for $Q^\s_{\leq
i}$, and $ \{ u_1^mu_2^j\tilde{\theta}^kv_1^l\, | \, m,j,k\geq 0,
\; l=0,1\}$ is an $\FF$-basis for $Q^\s$.
\end{lemma}

\Q
{\em Proof}.~~ This follows directly from the equality
\[
Q^\s = U +
Uv_1 = \bigoplus_{i,j\geq 0}\, \{ K[u_1] +
K[u_1]v_1\}u_2^i\tilde{\theta}^j
\]
and the fact that the leading terms of the elements $u_2,
\tilde{\theta}$ and $z_1$ are algebraically independent,  and from
the fact that $v_1 = z_1^3 +\cdots $ and $u_1 = z_1^2 +\cdots$.
{\smallbox}
%
%
\begin{theorem}\label{n=4has5gens}
The fixed algebra $F = K[x_1,x_2,x_3,x_4,x_5]^\s$
is generated by the five  elements
\[
x_1,\quad p_{1},\quad q_1,\quad p_2:= x_1^2u_2 = x_3^2 - x_2(x_3 + 2x_4)
+ x_1(x_3 + 3x_4 +2x_5)
\]
and
\begin{multline}
t:= x_1^3\tilde{\theta} = 6x_2^2x_5 - 6x_2(x_3 - x_2)x_4 + 2x_3^3 -
x_2x_3(3x_3 -x_2)\\ + x_1\left\{ -6(x_2 + 2x_3)x_5 + 9x_4^2 -6x_2x_4 +
2x_3^2 - x_2x_3\right\}
\end{multline}
of degrees $1,2,3, 2, 3$, respectively, subject to the  defining
relation
%
\begin{align}\label{relation1}
x_1^3t &= q_1^2 -p_1^3 + 3x_1p_1q_1 +2x_1^2p_1^2 + 3x_1^2p_1p_2.
\end{align}
For $i\geq 0$,
$$F_i = x_1^i Q_{\leq i}^\s =
\bigoplus \{ Kx_1^ap_1^bp_2^ct^dq_1^e\, | \, a+2b+2c+3t+3e=  i,
\, e=0,1\}$$ where $a,b,c,d\geq 0$.
\end{theorem}

\Q {\em Proof}.~~ By Corollary \ref{ngeq3} and Lemma \ref{f7}, $F
= \oplus_{i\geq 0}\, F_i$, where $$F_i = x_1^i Q_{\leq i}^\s =
\bigoplus \{ Kx_1^ap_1^bp_2^ct^dq_1^e\, | \, a+2b+2c+3t+3e=  i,
\, e=0,1\}.$$

By (\ref{tilde}), the relation (\ref{relation1}) holds. This
relation is {\em the} defining relation for the algebra $F$ since
the algebra that satisfies this defining relation has the same
basis as $F$ (see (\ref{relation1})). {\smallbox}

Finally, we make a comment on the case where $n\geq 5$.
Assume that $n\geq 5$.  Set $m:= \left[
\frac{n}{2}\right]$ and $ \mu:= \left[ \frac{n-1}{2}\right]$; then
$m =\mu$, when $n= 2m+1$, and $m=\mu +1$, when $n=2m$.   The
elements $u_k = u_k(z)$ and $v_k = v_k(z)$ of
Theorems~\ref{usubk},~\ref{tvsubk} can be written as sums $u_k =
u_k' + u_k''$ and $v_k = z_1u_k + v_k' +v_k''$, where $u_k', v_k'$
and $u_k'', v_k''$ are the {\em quadratic} terms and {\em linear}
terms, respectively. The elements

%
%
\begin{eqnarray}\label{wk}
w_k &:= v_k^2 - u_1u_k^2
    = u_k(2z_1v_k' -( z_1 + 2z_2)u_k) + (v_k')^2 + 2z_1u_kv_k'' +
2v_k'v_k'' + ( v_k'')^2,
\end{eqnarray}
for $k =2, \dots , \mu$, have degree $5$ and leading terms
%
%
\begin{eqnarray}\label{lwk}
l(w_k) = u_k'(2z_1v_k' - (z_1 +2z_2)u_k'),
\end{eqnarray}
for $k =2, \dots , \mu$.  Observe that $u_k'$ is the leading term of the
polynomial $u_k$.

\Q
The polynomials
%
%
\begin{eqnarray}\label{lastl}
u_k = (-1)^{k-1}2z_{2k} + \dots  \qquad {\rm and}\qquad v_k =
(-1)^{k-1}(1+2k)z_{2k+1} +\dots,
\end{eqnarray}
where by three dots we denote the terms from  $Q^{[2k-1]}$ and
$Q^{[2k]}$, respectively.  Hence,
%
%
\begin{eqnarray}\label{Q}
Q\equiv Q^{[n]} = Q^{[4]}[v_2, \dots,  v_{\mu}, u_3, \dots, u_m],
\end{eqnarray}
and
%
%
\begin{eqnarray}\label{Qsigma}
Q^{\sigma} = (Q^{[4]})^{\sigma} [v_2, \dots,  v_{\mu}, u_3, \dots,
u_m] = (Q^{[4]})^{\sigma}\otimes V,
\end{eqnarray}
where $V:= K[v_2, \dots  , v_\mu, u_3, \dots, u_m]$. Note that
\begin{eqnarray*}
(Q^{[5]})^{\sigma}&=& (Q^{[4]})^{\sigma} [v_2], \\
(Q^{[6]})^{\sigma}&=&(Q^{[4]})^{\sigma}[v_2, u_3],\\
(Q^{[7]})^{\sigma}&=&(Q^{[4]})^{\sigma}[v_2, v_3, u_3],\\
(Q^{[8]})^{\sigma}&=&(Q^{[4]})^{\sigma}[v_2, v_3, u_3, u_4].
\end{eqnarray*}
 The polynomial ring $V$ contains the polynomial subalgebra $W:=
K[w_2, \dots  , w_\mu, u_3, \dots, u_m]$ and
%
%
\begin{eqnarray}\label{V}
V= \bigoplus_{I\in \mZ_2^{\mu-1}}\, Wv^I,
\end{eqnarray}
is a free $W$-module of rank $2^{\mu-1}$, where $\mZ_2^{\mu-1} =
\mZ_2 \times \dots \times \mZ_2$ is the product of $\mu -1$ copies
of $\mZ_2 = \{0,1\}$, and $v^I = v_2^{i_2} \dots
v_{\mu}^{i_{\mu}}$, for $I= (i_3, \dots, i_\mu) \in \mZ_2^{\mu
-1}$.

\Q The leading terms of the polynomials $w_2,u_3, w_3, u_4, \dots,
$ are {\em algebraically independent}, since
%
%
\begin{eqnarray}\label{lastq}
l(u_k) = u_k' =  (-1)^{k-1}2z_1z_{2k-1} + \dots  \qquad {\rm
and}\qquad l(v_k')  = (-1)^{k-1}(2k-1)z_1z_{2k} +\dots,
\end{eqnarray}
where three dots denotes terms from $Q^{[2k-2]}$ and $Q^{[2k-1]}$,
respectively.
Thus, the sum
%
%
\begin{eqnarray}\label{W}
W = \sum Ku_3^{i_3} \dots u_m^{i_m}w_2^{j_2} \dots
w_{\mu}^{j_{\mu}},
\end{eqnarray}
for $i_2, \dots, j_\mu \geq 0$, is an $\FF$-direct sum.

\begin{lemma}
The sum (\ref{V}) is an $\FF$-direct sum in the filtered algebra $Q$ with
the filtration $\FF = \{ Q_{\leq i}\}$.
\end{lemma}

\Q
{\em Proof}.~~ This is evident, since $l(v_k) = z_1l(u_k), l(u_k) =
u_k'$ and $l(w_k)$ is as defined in (\ref{lwk}). {\smallbox}

\Q This means that one cannot produce new invariants from the
elements of $V$; that is, in any new invariants generators of
$Q^{[4]}$ must necessarily occur.


$${\bf Acknowledgements}$$

The authors would like to thank Claudio Procesi for giving them a
series of talks on invariant theory during his visit to the
University of Edinburgh in May 1999 and to Hanspeter Kraft for
explaining to the authors certain aspects of invariant theory in
the summer of 1998.




\bigskip


\noindent
V V Bavula \\
Department of Pure Mathematics\\
University of Sheffield\\
Hounsfield Road\\
Sheffield S3 7RH\\
Email: v.bavula@sheffield.ac.uk\\

\noindent T H Lenagan \\
School  of Mathematics\\
University of Edinburgh\\
James Clerk Maxwell Building\\
King's Buildings\\
Mayfield Road\\
Edinburgh EH9 3JZ\\
Email: tom@maths.ed.ac.uk


\end{document}